\documentclass[11pt,leqno]{amsart}
\usepackage{amssymb, amsmath}
\usepackage{pstricks,pst-node}
\usepackage{psfrag,graphicx}
\usepackage[all]{xy}

\let\oldsection=\section

\makeatletter
\renewcommand{\@seccntformat}[1]{\bf\@nameuse{the#1}.\quad}
\renewcommand\section{\@startsection{section}{1}%
            \z@{.7\linespacing\@plus\linespacing}{.5\linespacing}%
            {\normalfont\bfseries \boldmath}}
\renewcommand\subsection{\@startsection{subsection}{2}%
            \z@{.5\linespacing\@plus.7\linespacing}{-.5em}%
            {\normalfont\bfseries \boldmath}}
\renewcommand\subsubsection{\@startsection{subsubsection}{3}%
            \z@{.3\linespacing\@plus.5\linespacing}{-.5em}%
            {\normalfont\bfseries \boldmath}}
\makeatother

\theoremstyle{plain}
\newtheorem*{thm}{Theorem}
\newtheorem*{cor}{Corollary}
\newtheorem*{lem}{Lemma}
\newtheorem*{prop}{Proposition}

\theoremstyle{definition}

\newtheorem*{rems}{Remarks}

\numberwithin{equation}{subsection}

\newcounter{listequation}

\def\note#1{{\small\tt <<#1>>}}  
\def\note#1{}              

\def\Z{\mathbb Z}

\def\:{\colon}
\def\a{\alpha}
\def\b{\beta}

\def\la{\lambda}

\def\phi{\varphi}

\def\C{\mathbb C}
\def\CC{\mathbb C}

\def\:{\colon}

\def\JWw0S{{}^J W^{-w_0 S}}
\def\w0JWS{{}^{-w_0 J} W^S}
\def\fg{{\mathfrak g}}
\def\fh{{\mathfrak h}}
\def\fb{{\mathfrak b}}

\def\fn{{\mathfrak n}}

\def\rad{\operatorname{rad}}
\def\hd{\operatorname{hd}}
\def\soc{\text{soc}}
\def\Hom{\operatorname{Hom}}
\def\Ext{\operatorname{Ext}}
\def\End{\operatorname{End}}
\def\Im{\operatorname{Im}}

\def\sgn{\operatorname{sgn}}

\def\O{\mathcal O}
\def\W{\mathcal W}

\def\TJJ'{T_J^{J'}}
\def\TJ'J{T_{J'}^J}

\def\gr{\operatorname{gr}}
\def\CC{\mathcal C}
\def\op{\operatorname{op}}

\def\Label{\label}

\begin{document}
\title[Category $\O$ for the Virasoro algebra: Cohomology and Koszulity]
{Category $\O$ for the Virasoro algebra: Cohomology and Koszulity}

\author{Brian D. Boe}
\thanks{Research of the first author partially supported by NSA grant H98230-04-1-0103}
\address{Department of Mathematics \\
                      University of Georgia\\ Athens, Georgia 30602  }
\email{brian@math.uga.edu, nakano@math.uga.edu}

\author{Daniel K. Nakano}
\thanks{Research of the second author partially supported by NSF grant
DMS-0400548}

\author{Emilie Wiesner}
\address{Department of Mathematics \\
           Ithaca College\\ Ithaca, NY 14850  }
\email{ewiesner@ithaca.edu}
\thanks{Final version, to appear in Pacific Journal of Mathematics}
\date{\today}



\dedicatory{}

\begin{abstract}
In this paper the authors investigate blocks of the Category $\O$
for the Virasoro algebra over $\C$. We demonstrate that the blocks
have Kazhdan-Lusztig theories, and that the truncated blocks give rise to interesting 
Koszul algebras. The simple modules
have BGG resolutions, and from this the extensions between Verma modules and simple modules,
and between pairs of simple modules, are computed.
\end{abstract}

\maketitle

\parskip=2pt

\section{Introduction}
\Label{S:Intro}

In 1976, Bernstein, Gelfand and Gelfand \cite{BGG:76} initiated the study of
Category ${\mathcal O}$ for complex semisimple Lie algebras. Since the introduction
of Category ${\mathcal O}$, much progress has been made in studying the structure of
blocks for this category and its variants.
One of the major results in this area was the formulation and
proof of the Kazhdan-Lusztig (KL) Conjectures \cite{KaLu:79, BeBe:81, BrKa:81} which provide a
recursive formula for the characters of simple modules in Category ${\mathcal O}$.
These conjectures have been equivalently formulated in terms of $\text{Ext}$-vanishing
conditions between simple modules and Verma modules. For semisimple algebraic groups over fields of
positive characteristic $p>0$, an analogous conjecture has been provided by Lusztig
as long as  $p$ is at least as large as the Coxeter number of the underlying
root system. The characteristic $p$ Lusztig Conjecture still remains open.

In an attempt to better understand both the original Kazhdan-Lusztig
Conjecture and the Lusztig Conjecture, Cline, Parshall,
and Scott \cite{CPS1, CPS3, CPS5} developed an axiomatic treatment of highest weight categories
with the added structures involving ``Kazhdan-Lustzig theories'' and Koszulity.
Irving has partially developed some theories along these lines (see
\cite{Irv:90, Irv:92}). Beilinson, Ginzburg and Soergel \cite{BGS:96} proved that the principal
block of Category ${\mathcal O}$ is Koszul using perverse sheaves and
established Koszul duality between various blocks of Category ${\mathcal O}$, 
which provides an alternative proof of the KL Conjecture.
The work of Cline, Parshall and Scott is important because it isolates the key 
homological criteria for verifying the existence of such properties.

The Virasoro algebra is the universal central extension of the 
Witt algebra and plays a significant role in the definition 
of the vertex operator algebra. The theory of vertex operator 
algebras has provided a mathematical foundation for conformal 
field theory (cf.\ \cite{Lep}). Understanding such field theories in two 
dimensions involves problems about the representation theory of the 
Virasoro and vertex operator algebras. 

The Witt algebra is an infinite-dimensional simple Lie algebra over ${\mathbb C}$
and is the smallest example of a Cartan-type Lie algebra. The Virasoro algebra has a triangular decomposition $\fg=\fn^- \oplus \fh \oplus \fn^+$, which allows one to define a Category $\O$. In this
paper we study blocks of Category ${\mathcal O}$ for the Virasoro algebra, building on the foundational work of Feigin-Fuchs \cite{FeFu}, who determined all maps between Verma modules for the Virasoro algebra. After making explicit the construction of BGG-resolutions for simple modules in these blocks (which is implicit in \cite{FeFu}), we compute the ${\mathfrak n}^{+}$-cohomology with coefficients in any simple module. This extends results of Goncharova \cite{Gon1, Gon2}, who calculated
$\text{H}^{\bullet}(\fn^{+}, \C)$,
and of Rocha-Caridi and Wallach \cite{RCW2}, who computed  $\text{H}^{\bullet}(\fn^{+}, L)$
for any simple module
$L$ in the trivial block. This cohomological information allows us to calculate the extensions between
simple and Verma modules. We then verify that our categories satisfy properties given
in \cite{CPS5}; in particular, they have a KL theory. These properties yield a computation of 
extensions between all simple modules and imply that truncated blocks of Category ${\mathcal O}$ 
for the Virasoro algebra give rise to interesting Koszul algebras. We find it quite remarkable 
that KL theories naturally arise in the representation theory of the Virasoro algebra. It 
would be interesting to determine if this occurs in a more general context within the 
representation theory of Cartan-type Lie algebras.

The authors would like to thank Brian Parshall for conversations pertaining to calculating 
extensions in quotient categories, Jonathan Kujawa
for clarifying the connections between the extension theories 
used in Section \ref{S:Truncated}, and the referee for several helpful comments and suggestions.


\section{Notation and Preliminaries} \Label{S:Preliminaries}

\subsection{} The {\it Virasoro algebra} is the Lie algebra ${\mathfrak g} =
\C \mbox{-span} \{ z, d_k \mid k \in \Z \} $ with bracket $[\ ,\ ]$ given by
$$
[ d_k, z ] = 0, \quad [ d_j, d_k ] = (j-k)d_{j+k} + \frac{\delta_{j,-k}}{12}(j^3-j)z \qquad \text{for all } j, k\in \Z.
$$
The Virasoro algebra can be decomposed into a direct sum of subalgebras ${\mathfrak g} =
\fn^{-}\oplus \fh \oplus \fn^+ = \fn^- \oplus \fb^+$, where
$$\fn^-  =  \C\mbox{-span} \{d_n \mid n \in \Z_{<0} \}; \quad \fh = \C\mbox{-span} \{ d_0,z \}; \quad \fn^+ = \C\mbox{-span} \{d_n \mid n \in \Z_{>0} \},$$
and ${\mathfrak b}^{+}={\mathfrak h}\oplus {\mathfrak n}^{+}$. There is an anti-involution 
$\sigma: {\mathfrak g} \rightarrow {\mathfrak g}$ given by $\sigma(d_n)=d_{-n}$ and $\sigma(z)=z$.

\subsection{\bf Category $\mathcal O$ and Other Categories} \Label{SS:CategoryO}
The {\it Category $\mathcal{O}$} consists of ${\mathfrak g}$-modules $M$ such that
\begin{itemize}
\item $M= \bigoplus_{\mu \in \fh^{*}} M^{\mu}$, where $\fh^{*} = \Hom_{\C} (\fh, \C)$ and $M^{\mu} = \{m \in M \mid hm= \mu(h) m \ {\rm for \ all} \
h \in \fh \}$;
\item $M$ is finitely generated as a $\fg$-module;
\item $M$ is $\fn^+$-locally finite.
\end{itemize}
This definition is more restricted than the definition given in \cite{MoPi:95}.
Identify each
integer $n \in \Z$ with a weight $n \in \fh^{*}$ by $n(d_0) =n$ and $n(z)=0$.  Define a partial ordering on $\fh^{*}$ by 
\begin{equation} \Label{E:AltOrdering}
\mbox{$\mu < \gamma$ if  $\mu = \gamma + n$ for some $n \in \Z_{> 0}.$}
\end{equation}
The category defined in \cite{MoPi:95}, which we denote $\tilde{\O}$, consists of  $\fg$-modules $M$ such that $M=\bigoplus_{\mu \in \fh^{*}} M^{\mu}$, $\dim M^{\mu}<\infty$, and there exist $\lambda_1, \ldots \la_n \in \fh^*$ such that $M^{\mu} \neq 0$ only for $\mu \leq \lambda_i$ for some $i$.  Then $\O$ (as we have defined it) is the full subcategory of $\tilde{\O}$ consisting of finitely generated modules.  Therefore, many of the results about $\tilde{\O}$ proven in \cite{MoPi:95} apply to $\O$.  

For $\mu \in \fh^{*}$,  the {\it Verma module} $M(\mu)$ is the induced module
$$
M(\mu) = U({\mathfrak g}) \otimes_{U(\fb^+)} \C_{\mu}.
$$
The Verma module $M(\mu)$ has a unique simple quotient, denoted $L(\mu)$.  The modules $L(\mu)$,
$\mu \in \fh^{*}$, provide a complete set of simple modules in Category $\mathcal{O}$ (cf.\ \cite[Section 2.3]{MoPi:95}).
 For $\mu, \gamma \in \fh^{*}$, define a partial ordering 
 \begin{equation} \Label{E:Ordering}
 \mu \preceq \gamma \quad \mbox{if $L(\mu)$ is a subquotient of $M(\gamma)$}.
 \end{equation}  
 Extend this to an equivalence relation $\sim$.  The {\it blocks} of ${\mathfrak g}$ are the equivalence classes of $\fh^{*}$ determined by $\sim$.   For each block $[\mu] \subset \fh^{*}$, define $\mathcal{O}^{[\mu]}$ to be the full subcategory of $\mathcal{O}$ so that for $M \in \mathcal{O}^{[\mu]}$, where $L(\gamma)$ is a subquotient of $M$ only for $\gamma \in [\mu]$. For $M \in \mathcal{O}$, $M= \bigoplus_{[\mu] \subset \fh^{*}} M^{[\mu]}$, where $M^{[\mu]} \in \mathcal{O}^{[\mu]}$ (cf.\ \cite[2.12.4]{MoPi:95}). 

 We will make use of several other categories during the course of the paper, which we introduce now.  Let $\W$ be the category whose objects are ${\mathfrak g}$-modules $M$ such that $M= \bigoplus_{\lambda \in \fh^{*}} M^{\lambda}$; $M^{\lambda}$ is not necessarily finite-dimensional.   For a fixed weight $\mu \in \fh^*$, define $\W(\mu)$ to be the full subcategory of $\W$ whose objects are $\fg$-modules $M$ so that $M= \bigoplus_{\la \leq \mu} M^{\la}$.  

The anti-involution $\sigma$ can be used to define a duality functor $D$ on $\W$.
For $M \in \W$, define $D M = \bigoplus_\mu (M^{\mu})^{*}$ (as a vector space)
with ${\mathfrak g}$-action given by $(x f) (v)= f(\sigma(x)v)$ for $x \in {\mathfrak g}$,
$f \in DM$ and $v \in M$.  Then, $\Hom_\W (M,M') \cong \Hom_\W (DM', DM)$ for all $M, M' \in \W$.  Since $\sigma(h)=h$ for $h \in \fh$, $DM$ decomposes as a direct sum of weight spaces where $(DM)^{\mu} = (M^{\mu})^{*}$.  Therefore, $\W(\mu)$ is closed under  $D$.  Finally, note that $D L \cong L$ for any simple module $L \in \O$.

\section{ BGG Resolutions and $\fn^{+}$-Cohomology}
\Label{S:BGG}

\subsection{}  Theorem \ref{Virembedd} below, due to Feigin and Fuchs \cite{FeFu},  gives a description of
all Verma module embeddings in a given block of $\mathcal O$.  Since every non-zero map between Verma modules is an embedding, this describes all homomorphisms between Verma modules in a block. From this result one can construct
BGG resolutions of the simple modules $L(\mu)$ to compute $\text{H}^{\bullet}(\fn^+, L(\mu))$.

\begin{thm}[{\cite[1.9]{FeFu}}] \Label{Virembedd} Suppose $\mu \in \fh^{*}$, and set $h=\mu(d_0), c=\mu(z) \in \C$.  Define
$$
\nu=\frac{c-13+ \sqrt{(c-1)(c-25)}}{12}; \quad \beta= \sqrt{-4 \nu h + (\nu+1)^2}
$$
and consider the line in the $rs$-plane
\begin{equation}
\mathcal{L}_{\mu}: r+ \nu s + \beta=0.
\end{equation}
  The Verma module embeddings involving $M(\mu)$ are determined by integer points $(r,s)$ on $\mathcal{L}_{\mu}$:
\begin{itemize}
\item[(i)] Suppose $\mathcal{L}_{\mu}$ passes through no integer points or one integer point $(r,s)$ with $rs=0$.  Then the block $[\mu]$ is given by $[\mu]= \{ \mu \}.$

\item[(ii)] Suppose $\mathcal{L}_{\mu}$ passes through exactly one integer point $(r,s)$ with $rs \neq 0$.  The block $[\mu]$ is given by $[\mu] = \{ \mu, \mu +rs \}.$  The block structure is given by Figure \ref{Crslittle}, where an arrow $\la \rightarrow \gamma$ between weights indicates $M(\la) \subseteq M(\gamma)$.

\psfrag{Mmu}{{\footnotesize{$\mu=\mu_0$}}}
\psfrag{mu+rs,1}{{\footnotesize $\mu_1=\mu+rs$}}
\psfrag{mu+rs,-1}{{\footnotesize $\mu_{-1}=\mu+rs$}}
\psfrag{rs>0}{{\footnotesize $rs>0$}}
\psfrag{rs<0}{{\footnotesize $rs<0$}}
\psfrag{-1}{{\footnotesize $(r_{-1}, s_{-1})$}}
\psfrag{1}[Br]{{\footnotesize $(r_{1}, s_{1})$}}
\psfrag{2}{{\footnotesize $(r_{2}, s_{2})$}}
\psfrag{3}{{\footnotesize $(r_{3}, s_{3})$}}
\psfrag{4}[Br]{{\footnotesize $(r_{4}, s_{4})$}}
\psfrag{5}[Br]{{\footnotesize $(r_{5}, s_{5})$}}
\psfrag{6}{{\footnotesize $(r_{6}, s_{6})$}}
\psfrag{7}{{\footnotesize $(r_{7}, s_{7})$}}
\psfrag{8}[Br]{{\footnotesize $(r_{8}, s_{8})$}}
\psfrag{9}[Br]{{\footnotesize $(r_{9}, s_{9})$}}
\psfrag{slope1}{{\footnotesize$slope \left( \mathcal{L}_{\mu} \right)>0$}}
\psfrag{slope2}{{\footnotesize$slope \left( \mathcal{L}_{\mu} \right)<0$}}
\begin{figure}[h]
\centering
\includegraphics[height=.08\textheight]{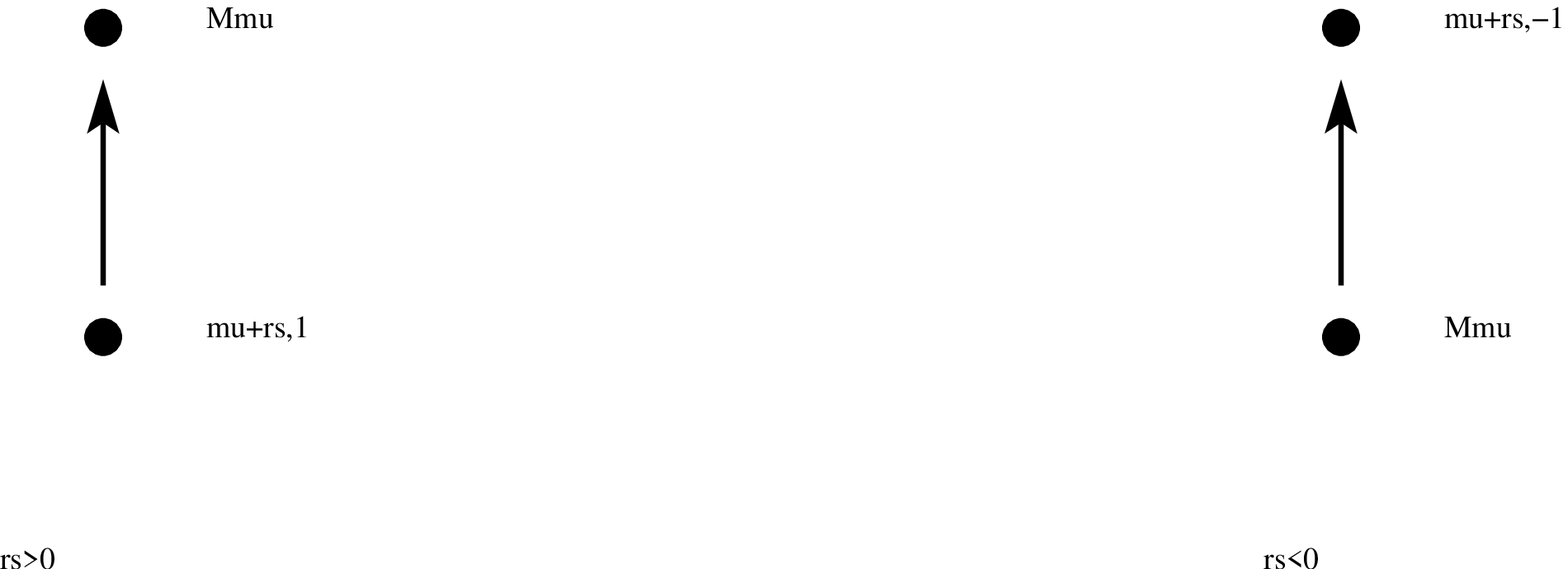}
\caption{Case (ii)} \Label{Crslittle}
\end{figure}

\item[(iii)] Suppose $\mathcal{L}_{\mu}$ passes through infinitely many integer points and crosses an axis at an integer point.  Label these points $(r_i,s_i)$ so that 
$$\cdots < r_{-2}s_{-2}<r_{-1}s_{-1}<0=r_0s_0<r_1s_1<r_2s_2 <\cdots$$
as in Figure \ref{Lmu}.   The block $[\mu]$ is given by $ [\mu] = \{ \mu_i=\mu+r_is_i  \}$.  The block structure is given by  Figure \ref{Crs1}.

\psfrag{slope1}{{\footnotesize$slope \left( \mathcal{L}_{\mu} \right)>0$}}
\psfrag{slope2}{{\footnotesize$slope \left( \mathcal{L}_{\mu} \right)<0$}}
\psfrag{VM12}{{\footnotesize $\mu=\mu_0$}}
\psfrag{VM30}{{\footnotesize $\mu=\mu_0$}}
\psfrag{VM4}{{\footnotesize $\mu=\mu_0$}}
\psfrag{VM5}{{\footnotesize ${\mu_{-1}}$}}
\psfrag{VM6}{{\footnotesize$\mu_1$}}
\psfrag{VM15}[Br]{{\footnotesize$\mu_{-2}$}}
\psfrag{VM16}[Br]{{\footnotesize$\mu_{-1}$}}
\psfrag{VM18}[Br]{{\footnotesize$\mu_1$}}
\psfrag{VM19}[Br]{{\footnotesize$\mu_2$}}
\psfrag{0}[Br]{{\footnotesize $(r_{0}, s_{0})$}}
\begin{figure}
\begin{minipage}[b]{0.45\linewidth}
\centering
\includegraphics[height=.175\textheight]{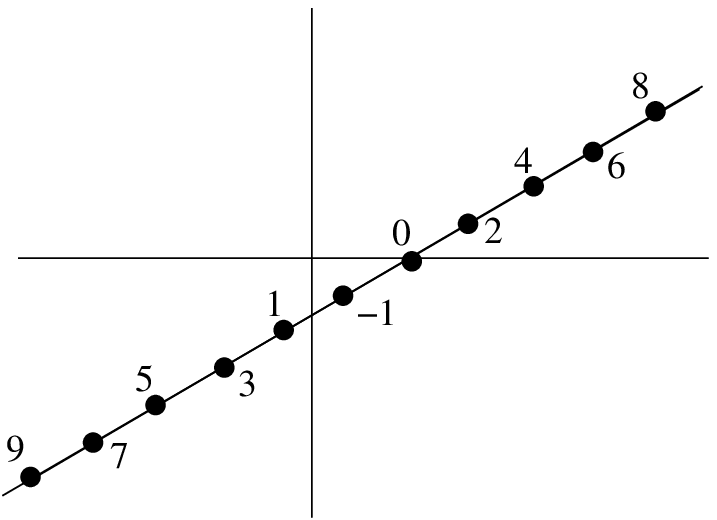}
\caption{$\mathcal{L}_{\mu}$} \Label{Lmu}
\end{minipage}
\begin{minipage}[b]{0.45\linewidth} 
\centering
\includegraphics[height=.175\textheight]{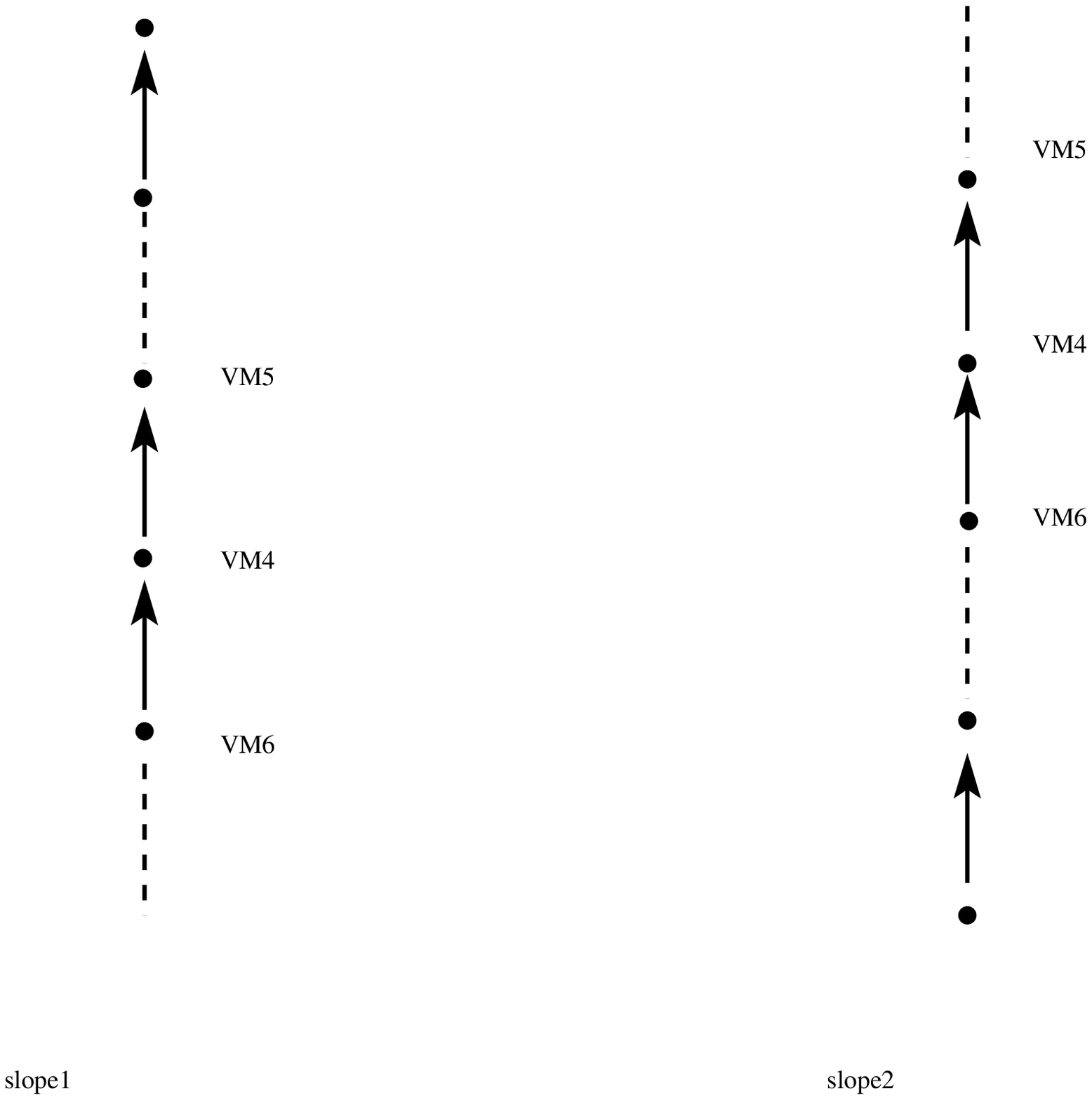}
\caption{Case (iii)} \Label{Crs1}
\end{minipage}
\end{figure}

\item[(iv)] Suppose $\mathcal{L}_{\mu}$ passes through infinitely many integer points and does not cross either axis at an integer point.  Label the integer points $(r_i,s_i)$ on $\mathcal{L}_{\mu}$ so that
$$\cdots <r_{-1}s_{-1}<r_0s_0<0<r_1s_1<r_2s_2 <\cdots.$$   
Also consider the auxiliary line $\tilde{\mathcal{L}}_{\mu}$ with the same slope as $\mathcal{L}_{\mu}$ passing through the point $(-r_1,s_1)$.  Label the integer points on this line $(r_j',s_j')$ in the same way as $\mathcal{L}_{\mu}$.
The block $[\mu]$ is given by $[\mu]=\{ \mu_i, \mu_i' \}$, where
$$
\mu_i = \left\{
\begin{array}{ll}
\mu+r_is_i, & i \ \mbox{odd} \\
\mu+r_1s_1+ r_{i}'s_{i}', & i \ \mbox{even}
\end{array},
\right.
\quad
\mu_i' = \left\{
\begin{array}{ll}
\mu+r_{i+1}s_{i+1}, & i \ \mbox{odd} \\
\mu+r_1s_1+ r_{i+1}'s_{i+1}', & i \ \mbox{even}
\end{array}.
\right.
$$
The block structure is given by Figure \ref{embedd4}.

\psfrag{slope1}{{\footnotesize$slope \left( \mathcal{L}_{\mu} \right)>0$}}
\psfrag{slope2}{{\footnotesize$slope \left( \mathcal{L}_{\mu} \right)<0$}}
\psfrag{VM12}{{\footnotesize $\mu=\mu_0$}}
\psfrag{VM30}{{\footnotesize $\mu=\mu_0$}}
\psfrag{VM4}{{\footnotesize $\mu=\mu_0$}}
\psfrag{VM5}{{\footnotesize ${\mu_{-1}'}$}}
\psfrag{VM6}{{\footnotesize$\mu_1'$}}
\psfrag{VM10}{{\footnotesize$\mu_{-2}$}}
\psfrag{VM11}{{\footnotesize$\mu_{-1}$}}
\psfrag{VM13}{{\footnotesize$\mu_1$}}
\psfrag{VM14}{{\footnotesize$\mu_2$}}
\psfrag{VM15}[Br]{{\footnotesize$\mu_{-2}'$}}
\psfrag{VM16}[Br]{{\footnotesize$\mu_{-1}'$}}
\psfrag{VM17}[Br]{{\footnotesize$\mu_0'$}}
\psfrag{VM18}[Br]{{\footnotesize$\mu_1'$}}
\psfrag{VM19}[Br]{{\footnotesize$\mu_2'$}}
\begin{figure}[h]
\centering
\includegraphics[height=.24\textheight]{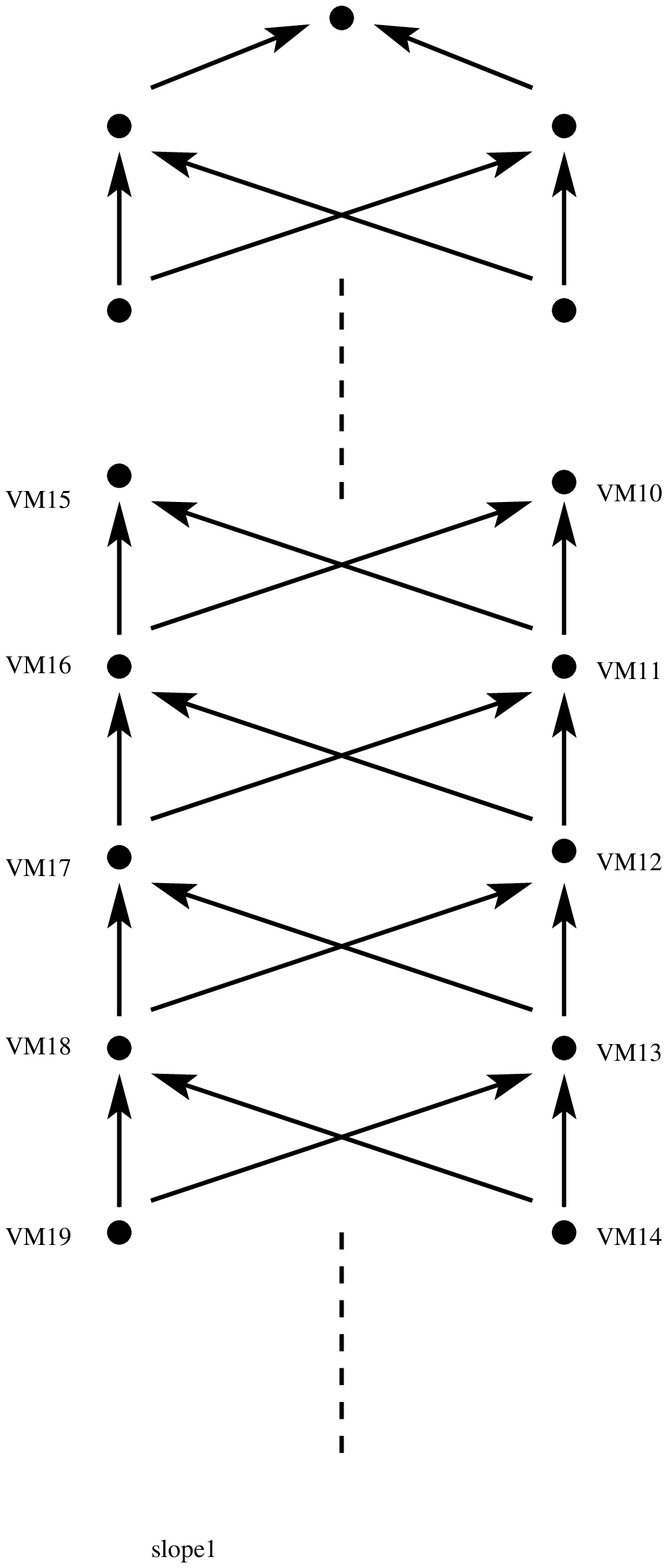}
\hspace{2in}
\includegraphics[height=.24\textheight]{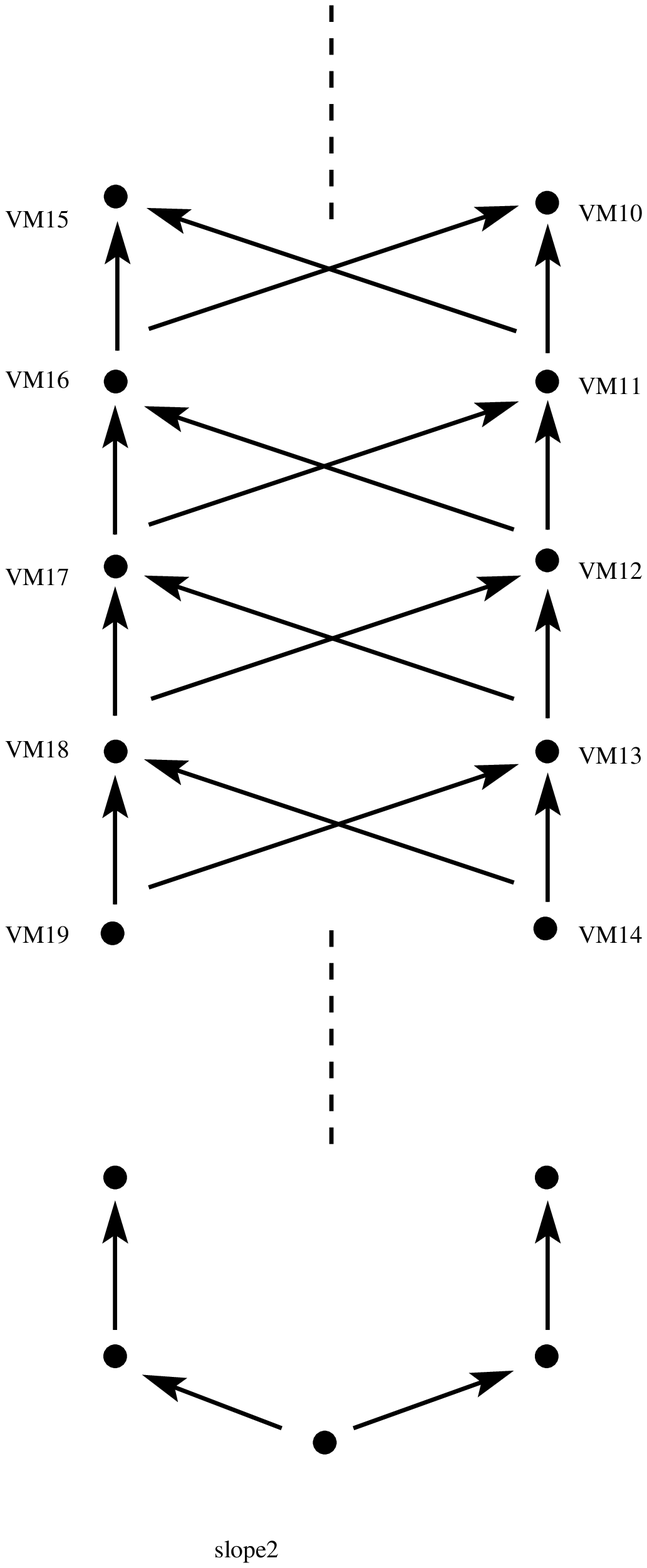}
\caption{Case (iv)} \Label{embedd4}
\end{figure}
\end{itemize}
\end{thm}
We will refer to blocks as in case (iii) as ``thin'' blocks and blocks as in case (iv) as ``thick'' blocks.
The second type of thick block has a highest weight poset structure equivalent to the Bruhat order on $D_\infty$, the infinite dihedral group.

\subsection{\bf BGG Resolutions} Feigin and Fuchs \cite{FeFu2}  observe without elaboration that their result (Theorem \ref{Virembedd}, above) yields a BGG resolution for the simple modules in Category $\mathcal{O}$.  We now provide the details for constructing these resolutions.

Given a module $M$ in Category $\O$, define the \emph{radical} of $M$, $\rad M$, to be the smallest submodule such that $M/\rad M$ is semisimple. Put $\rad^0 M=M$ and for $i>0$, $\rad^i M=\rad(\rad^{i-1} M)$. This defines a decreasing filtration of $M$, the \emph{radical filtration}. For $i\ge 0$, \emph{layer} $i$ of the radical filtration is defined to be $\rad_i M = \rad^i M/\rad^{i+1} M$. We also write $\hd M = M/\rad M$.
In general (cf.\ \cite{FeFu}), the terms and the layers of the radical filtration of $M(\mu)$, in the notation of Figures \ref{Crslittle}, \ref{Crs1}, and \ref{embedd4}, are as follows:
\begin{equation} \Label{E:VMrad}
\begin{gathered}
\rad^i M(\mu) = M(\mu_i) + M(\mu'_i) \ \text{ for } i>0,\cr
\rad_0 M(\mu) = L(\mu), \quad \rad_i M(\mu) = L(\mu_i) \oplus L(\mu_i') \ \mbox{ for} \ i>0.
\end{gathered}
\end{equation}
If $\mu$ belongs to a finite or thin block, then terms involving $\mu'_i$ are to be ignored. Also if $\mu$ belongs to a block with a minimal element, say $\mu_n$, then 
$$
\rad^n M(\mu) = \rad_n M(\mu) = M(\mu_n) = L(\mu_n), \quad \rad^i M(\mu)=\rad_i M(\mu) = 0 \ \mbox{for} \ i>n.
$$

Assume $\mu$ belongs to a thick block. According to \cite{BGG:75}, there will be a complex $C_\bullet \to L(\mu)\to 0$, where $C_i$ is the direct sum of the Verma modules $M(\mu_i)\oplus M(\mu'_i)$, provided that to each edge of the poset below $\mu$, using the ordering $\preceq$, it is possible to assign  a sign $+1$ or $-1$ in such a way that the product of the signs on any diamond is $-1$. Such a labeling is indicated in Figure \ref{F:thick}. It is easy to check directly in this case that the resulting complex is in fact a resolution, called a BGG resolution of $L(\mu)$.

\begin{figure}[ht]
\centering
\begin{pspicture}(-.3,-.35)(.2,7)
\psset{linewidth=.5pt,labelsep=8pt,nodesep=0pt}
\small
$
\rput(0,0){\rnode{a0}{\mu_n}}
\rput(-1,1){\rnode{a1}{\mu_{n-1}'}}
\rput(1,1){\rnode{b1}{\mu_{n-1}}}
\rput(-1,2){\rnode{a2}{\mu_{n-2}'}}
\rput(1,2){\rnode{b2}{\mu_{n-2}}}
\pnode(-1,2.5){p1}
\pnode(1,2.5){p2}
\pnode(-1,3){q1}
\pnode(1,3){q2}
\rput(-1,3.5){\rnode{a3}{\mu_{3}'}}
\rput(1,3.5){\rnode{b3}{\mu_{3}}}
\rput(-1,4.5){\rnode{a4}{\mu_2'}}
\rput(1,4.5){\rnode{b4}{\mu_2}}
\rput(-1,5.5){\rnode{a5}{\mu_1'}}
\rput(1,5.5){\rnode{b5}{\mu_1}}
\rput(0,6.5){\rnode{a6}{\mu}}
\everypsbox{\scriptstyle}
{\psset{linewidth=.5pt,nodesep=3pt,labelsep=2pt}
\ncline{a0}{a1}\Aput{+}
\ncline{a0}{b1}\Bput{\mp}
\ncline{a1}{a2}\Aput{+}
\ncline{a1}{b2}\bput(.22){\pm}
\ncline{b1}{a2}\aput(.2){\pm}
\ncline{b1}{b2}\Bput{+}
\ncline{a3}{a4}\Aput{+}
\ncline{a3}{b4}\bput(.22){+}
\ncline{b3}{a4}\aput(.18){+}
\ncline{b3}{b4}\Bput{+}
\ncline{a4}{a5}\Aput{+}
\ncline{a4}{b5}\bput(.22){-}
\ncline{b4}{a5}\aput(.18){-}
\ncline{b4}{b5}\Bput{+}
\ncline{a5}{a6}\Aput{+}
\ncline{b5}{a6}\Bput{+}
}
{\psset{linewidth=2pt,linestyle=dotted,dotsep=6pt}
\ncline{p1}{q1}
\ncline{p2}{q2}
}
$
\end{pspicture}
\hspace{1.5in}
\begin{pspicture}(-.3,-.35)(.2,7)
\psset{linewidth=.5pt,labelsep=8pt,nodesep=0pt}
\small
$
\pnode(0,2){p1}
\pnode(0,3){q1}
\rput(-1,3.5){\rnode{a3}{\mu_{3}'}}
\rput(1,3.5){\rnode{b3}{\mu_{3}}}
\rput(-1,4.5){\rnode{a4}{\mu_2'}}
\rput(1,4.5){\rnode{b4}{\mu_2}}
\rput(-1,5.5){\rnode{a5}{\mu_1'}}
\rput(1,5.5){\rnode{b5}{\mu_1}}
\rput(0,6.5){\rnode{a6}{\mu}}
\everypsbox{\scriptstyle}
{\psset{linewidth=.5pt,nodesep=3pt,labelsep=2pt}
\ncline{a3}{a4}\Aput{+}
\ncline{a3}{b4}\bput(.22){+}
\ncline{b3}{a4}\aput(.18){+}
\ncline{b3}{b4}\Bput{+}
\ncline{a4}{a5}\Aput{+}
\ncline{a4}{b5}\bput(.22){-}
\ncline{b4}{a5}\aput(.18){-}
\ncline{b4}{b5}\Bput{+}
\ncline{a5}{a6}\Aput{+}
\ncline{b5}{a6}\Bput{+}
}
{\psset{linewidth=2pt,linestyle=dotted,dotsep=6pt}
\ncline{p1}{q1}
}
$
\end{pspicture}

\caption{Assignment of signs}
\Label{F:thick}
\end{figure}
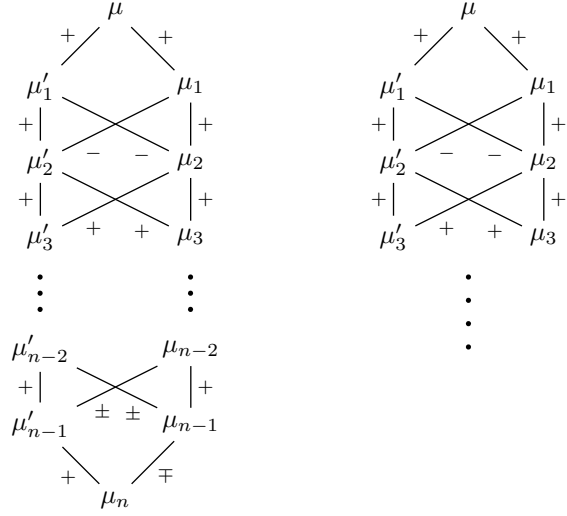

Explicitly, the following are BGG resolutions of  $L(\mu)$: for $\mu$ belonging to a thick block with a minimal element $\mu_{n}$,
$$
0 \to M(\mu_n) \to M(\mu_{n-1}) \oplus M(\mu_{n-1}') \to \cdots \to M(\mu_1) \oplus M(\mu_1') \to  M(\mu)  \to L(\mu) \to 0;
$$
and for $\mu$ belonging to a thick block with a maximal element,
$$
\cdots \to M(\mu_i) \oplus M(\mu_i') \to \cdots \to M(\mu_1) \oplus M(\mu_1') \to  M(\mu)  \to L(\mu) \to 0
$$

Next consider a weight $\mu \in \fh^{*}$ belonging to a  thin block or a finite block.  Then $\rad M(\mu) = M(\mu_1)$ if $\mu_1$ exists in the block, and $\rad M(\mu) = 0$ otherwise.  Thus, if $M(\mu)$ is not itself irreducible, the BGG resolution of $L(\mu)$ is
$$
0 \to M(\mu_1) \to M(\mu) \to L(\mu) \to 0.
$$

We now introduce some additional notation.  Fix $\mu \in \fh^{*}$.  Define a length function $l: [\mu] \rightarrow \Z$ by 
$$
l(\mu_i)=l(\mu'_i)=i
$$
(using the notation of Figures \ref{Crslittle}, \ref{Crs1}, and \ref{embedd4}).
While $l(\ )$ depends on a choice of representative $\mu$ for the block, the value (and, in particular, the parity) of $l(\nu)-l(\gamma)$  ($\nu, \gamma \in [\mu]$) is independent of the choice of representative.  This will be relevant later in the paper.

In summary of the description given above, each simple module $L(\mu)$ has a BGG resolution $\cdots \to C_1\to C_0 \to L(\mu) \to 0$ where 
\begin{equation} \Label{E:BGGres}
C_i = \left\{
\begin{array}{ll}
\  \bigoplus_{l(\nu)=i,\ \nu\preceq\mu} M(\nu) & \mbox{if $[\mu]$ is a thick block or a finite block}; \\[5pt]
\begin{array}{ll}
 \bigoplus_{l(\nu)=i} M(\nu), &i \leq 1\\
0, & i>1 
\end{array}
& 
 \mbox{if $[\mu]$ is a thin block}.
\end{array}
\right.
\end{equation}

\subsection{\bf $\fn^+$-Cohomology}
In  \cite{Gon1, Gon2}, Goncharova proved that 
\begin{equation}
\text{H}^k (\fn^{+}, \C) = \C_{\frac{3k^2+k}{2}}
\oplus \C_{\frac{3k^2-k}{2}}.
\end{equation}
Rocha-Caridi and Wallach  use Goncharova's work to obtain BGG resolutions
for $\C$ \cite{RCW2} and for the other simple modules in the trivial block \cite{RCW3}.  Using this,
they compute $\text{H}^k (\fn^{+}, L)$ for any simple module $L$ in the trivial block \cite{RCW3}.
We extend Rocha-Caridi and Wallach's result to cohomology with coefficients in any simple module in $\mathcal O$.

\begin{thm} \Label{T:n-homology}
Let $\mu \in \fh^{*}$, and let $k \in \Z_{\geq 0}$.  
\begin{itemize}
\item[(a)] \Label{E:n-hom-thick} Suppose that $\mu$ belongs to a thick block or a finite block.  As an $\fh$-module,
 $$
 \operatorname{H}^k(\fn^+,L(\mu))  \cong \operatorname{H}_k(\fn^-,L(\mu)) \cong \bigoplus_{\mbox{\scriptsize{\shortstack{$\nu \in [\mu]$, $\nu \preceq \mu$ \\ $l(\nu)-l(\mu)=k$}}}} \C_{\nu}.
 $$
\item[(b)]  \Label{E:n-hom-thin}
 Suppose that $\mu$ belongs to a thin block.  As an $\fh$-module,
$$
\operatorname{H}^k(\fn^+,L(\mu))  \cong \operatorname{H}_k(\fn^-,L(\mu)) \cong 
\left\{
\begin{array}{cl}
\displaystyle\bigoplus_{\mbox{\scriptsize{\shortstack{$\nu \in [\mu]$, $\nu \preceq \mu$ \\ $l(\nu)-l(\mu)=k$}}}}  \C_{\nu} & \mbox{if $k \leq 1$}; \\[30pt]
0 & \mbox{if $k>1$.}
\end{array}
\right. 
$$
\end{itemize}
\end{thm}
\begin{proof}
We first compute the homology groups $\text{H}_k(\fn^-,L(\mu))$, where 
$\text{H}_k(\fn^-,-)$ is the $k$th left derived functor
of $\C \otimes_{U(\fn^-)} -$.  Since Verma modules are free $U(\fn^-)$-modules,  apply
$\C \otimes_{U(\fn^-)} -$ to the resolution \eqref{E:BGGres}.  Note that $\C  \otimes_{U(\fn^-)} M(\nu) \cong \C_{\nu}$ .  The resulting differential maps in the resolution are 
${\mathfrak h}$-equivariant and $\C_{\nu}$ appears at most once in the resolution for each weight $\nu$.  Therefore, all of the differential maps must be zero.
This verifies that $\text{H}_k(\fn^-,L(\mu))$ is given by the formula 
in the statement of the theorem.

We now show that $\text{H}_k(\fn^-,L(\mu)) \cong \text{H}^k(\fn^+, L(\mu))$.  This may be well known; it is claimed in \cite{RCW2} to follow from ``standard arguments.'' Because the infinite dimensional case seems somewhat subtle, we include a proof for completeness. 

Write $L=L(\mu)$.  Recall that $\text{H}_k(\fn^-,L)$ can be computed using the complex
$$\cdots \rightarrow \Lambda^{k}(\fn^-) \otimes L \stackrel{d_k}{\rightarrow}
\Lambda^{k-1} (\fn^-) \otimes L \rightarrow \cdots$$
and $\text{H}^k(\fn^+,L)$ can be computed using the complex
$$\cdots \rightarrow \Lambda^k((\fn^+)^{*}) \otimes L \stackrel{d^k}{\rightarrow}
\Lambda^{k+1}((\fn^+)^{*}) \otimes L \rightarrow \cdots.$$

We extend the notion of duality defined in Section \ref{SS:CategoryO} to
$\Lambda^k(\fn^-)$, viewed as an $\fh$-module, as follows. 
For $f \in \oplus_{\la \in \fh^*} ((\Lambda^k(\fn^-))^{*})^{\la}$, define
$\tilde{f} \in ( \Lambda ^k (\fn^+))^*$ by $\tilde{f} (x) = f(\sigma(x))$ for $x \in \Lambda^k(\fn^+)$.  Let $D(\Lambda^k(\fn^-)) = \{ \tilde{f} \mid f \in  \oplus_{\la \in \fh^*} ((\Lambda^k(\fn^-))^{*})^{\la} \}$.  Define $\Lambda^k(D(\fn^-)) \subseteq \Lambda^k((\fn^+)^{*})$ analogously. 

By choosing a basis for each weight space $\left( \Lambda^k (\fn^-) \otimes L \right)^{\la}$, we can construct an $\fh$-module embedding 
$\Lambda^k (\fn^-) \otimes L \rightarrow D(\Lambda^k (\fn^-)) 
\otimes L \subseteq (\Lambda^k (\fn^+))^* \otimes L$. Since $D(L) \cong L$, the differential map $d_{k+1}$ induces a codifferential map $\tilde{d}_k:
D(\Lambda^k(\fn^-))  \otimes L \rightarrow D(\Lambda^{k+1}(\fn^-)) \otimes L$ as follows.  Let $f \in \oplus_{\lambda \in \fh^*} (\Lambda^k(\fn^-))^*)^{\lambda}$ and $g \in \oplus_{\lambda \in \fh^*} (L^*)^{\la}$.  Then $f \otimes g$ corresponds to an element $\widetilde{f \otimes g} \in D(\Lambda^k(\fn^-))  \otimes D L \cong D(\Lambda^k(\fn^-))  \otimes L$.   Define $\tilde{d}_k (\widetilde{f \otimes g})(x \otimes m) = (f \otimes g) (d_{k+1} ( \sigma(x) \otimes m))$ for $x \in \Lambda^{k+1} (\fn^+)$ and $m \in L$.  For $\lambda \in \fh^*$, it can be shown that $\dim (\ker \tilde{d}_k)^{\lambda} - \dim (\Im \tilde{d}_{k-1})^{\lambda}= \dim (\ker {d}_{k})^{\lambda} - \dim (\Im {d}_{k+1})^{\lambda}$.  This implies
\begin{equation} \label{E:n-cohom-proof1}
\text{H}^k(D(\Lambda^{\bullet}(\fn^-))\otimes L) \cong \text{H}_k(\fn^-,L).
\end{equation}

For each $k$ define $\phi_k=\phi:\Lambda^k((\fn^+)^{*}) \rightarrow \left( \Lambda^k(\fn^+) \right)^{*}$ by
$$
\phi(f_1 \wedge \cdots\wedge f_k) (x_1 \wedge \cdots \wedge x_k) =
\sum_{\tau \in S_k} \sgn(\tau) \prod_i f_i (x_{\tau(i)}), \quad\text{for } x_i \in \fn^+,\ f_i \in (\fn^+)^{*}.
$$
The map $\phi$ is an $\fh$-module isomorphism, and
$\phi(\Lambda^k(D(\fn^-))) = D(\Lambda^k(\fn^-))$.  Moreover, it can be checked that
$(\phi_{k+1}\otimes 1) \circ d^k = \tilde{d}_k \circ (\phi_k \otimes 1)$ on $D(\Lambda^k (\fn^-))\otimes L$.   Therefore, $\phi \otimes 1$ gives an isomorphism 
\begin{equation} \label{E:n-cohom-proof2}
\text{H}^k (D(\Lambda^{\bullet}(\fn^-))\otimes L) \cong \text{H}^k (\Lambda^{\bullet}(D (\fn^-))\otimes L) .
\end{equation}  
To complete the proof, we need to show that $\text{H}^k (\Lambda^{\bullet}(D (\fn^-))\otimes L) \cong \text{H}^k (\fn^+, L)$, which entails checking the following:
\begin{itemize}
\item[(i)] For $X \in \Im (d^{k-1}) \cap \Lambda^k (D (\fn^-))\otimes L$, there is $Y \in \Lambda^{k-1}(D (\fn^-))\otimes L$ with $d^{k-1}(Y)=X$.  
\item[(ii)] For $X \in \ker (d^k)$, there is $\tilde{X} \in \Lambda^k (D (\fn^-)) \otimes L$
such that $X - \tilde{X} \in \Im(d^{k-1})$.   
\end{itemize}

Let $X \in  \Lambda^k((\fn^+)^{*}) \otimes L$.  The space
$\Lambda^k((\fn^+)^{*}) \otimes L$ decomposes as a direct product of its weight spaces.
Write $X = \prod_{\lambda \in \fh^{*}} X_{\lambda}$, where
$X_{\lambda} \in (\Lambda^k((\fn^+)^{*}) \otimes L)^{\lambda}$.   Then
$X_{\lambda} = f_1 \otimes m_1 + \cdots +f_n \otimes m_n$ for some
$f_i \in \Lambda^k((\fn^{+})^{*})^{\mu_i}$ and $m_i \in L^{\nu_i}$ with 
$\mu_i + \nu_i = \lambda$. From the definition of $\Lambda^k(D(\fn^-))$ we see that $\Lambda^k(D(\fn^-))$ is the set of
elements $f \in \Lambda ^k ((\fn^{+})^{*})$, $f = \prod_{\mu \in \fh^{*}} f_{\mu}$, such that
$f_{\mu} \neq 0$ for only finitely many $\mu$.  Therefore, $f_i \in \Lambda^k
(D(\fn^-))$, and so
$X_{\lambda} \in \Lambda^k( D(\fn^-)) \otimes L$.

The differential map $d^k$ preserves weight spaces.  Suppose $X = d^{k-1}(\tilde{Y})$ for some $\tilde{Y} \in \Lambda^{k-1}((\fn^+)^{*}) \otimes L$.  Define $Y= \prod_{\lambda \in \fh^*} Y_{\lambda}$ by $Y_{\lambda}=\tilde{Y}_{\lambda}$ if $X_{\lambda} \neq 0$ and $Y_{\lambda}=0$ otherwise.  If $X \in \Lambda^k( D(\fn^-)) \otimes L$ then $Y \in \Lambda^{k-1}( D(\fn^-)) \otimes L$.  This proves (i).  

Now let $X \in \ker(d^k)$.   Then $d^k (X_{\la})=0$ for all $\la \in \fh^*$.  Equations (\ref{E:n-cohom-proof1}) and (\ref{E:n-cohom-proof2}) imply that $\text{H}^k (\Lambda^{\bullet}(D (\fn^-))\otimes L)$ is finite-dimensional.  Because  $X_{\lambda} \in \Lambda^k( D(\fn^-)) \otimes L$, this means that there are only finitely many  $\lambda_1,
\ldots,  \lambda_r$ such that $X_{\lambda_i} \neq 0$ and
$X_{\lambda_i} \neq d^k ( Y_{\lambda_i} )$ for some
$Y_{\lambda_i} \in \Lambda^{k-1}(D(\fn^-)) \otimes L$.  Define $\tilde{X}_{\lambda_i} =
X_{\lambda_i}$ for $i=1, \ldots, r$ and $\tilde{X}_{\lambda}=0$ otherwise.  Then
$\tilde{X} \in (\Lambda^k(D(\fn^-)) \otimes L$ and $X - \tilde{X} \in \Im (d^{k-1})$.  This proves (ii).
\end{proof}

\begin{cor}
Every irreducible module in Category $\O$ for the Virasoro algebra is a Kostant module (in the sense of \cite{BoHu}).
\end{cor}

\section{Extensions}

The structure of the infinite blocks of $\O$ presents various obstacles in computing $\Ext$-groups. The infinite blocks with a minimal element do not have enough projectives.  In the infinite blocks with a maximal element, objects do not generally have finite length.  We demonstrate that the first problem can be remedied by truncation, and the second can be addressed via a quotient construction.

\subsection{\bf Cohomology and Truncated Categories} \Label{S:Truncated} We first define the truncation of an infinite block of $\O$ having a minimal element.
Fix a weight $\mu$ in the block $\mathcal C$, and consider the full subcategory ${\mathcal C}(\mu)$---called the \emph{truncation of the block at $\mu$}---of modules all of whose composition factors have highest weights less than or equal to $\mu$, using the partial ordering given in (\ref{E:Ordering}).   

Now let $\mathcal C$ be a finite block of $\mathcal O$, an infinite block of $\O$ with a maximal element, or a truncated infinite block with minimal element.  Denote the weight poset of $\mathcal C$ by $\Lambda$.  Then there is a maximal element $\mu \in \Lambda$.  If $\mathcal C$ is a truncated thick block, we can write $\mu= \mu_0$ as in Figure \ref{embedd4}.  We assume that $\mu$ is chosen so that $\mu_0 \leq \mu_0'$ in the partial ordering given by \eqref{E:AltOrdering}.  Then $\Lambda = \{ \nu \in [\mu] \mid \nu \leq \mu \}$, which allows us to compare $\mathcal C$ and $\W(\mu)$.  
We now use Theorem \ref{T:n-homology} to compute $\Ext^i_{\mathcal C} (M(\lambda),L(\nu))$ by passing through relative cohomology and using the categories $\W$ and $\W(\mu)$ defined in Section \ref{SS:CategoryO}.
\begin{lem}[A] \Label{L:relativecohom}
Let $\lambda, \nu \in \Lambda$ and $i \in \Z_{\geq 0}$.  Then 
$\Ext_{\W}^i \left( M(\lambda), L(\nu) \right)\cong \operatorname{H}^i (\fn^+, L(\nu) )^{\lambda}$.
\end{lem}
\begin{proof}
For $i \in \Z_{\geq 0}$, define $P_i = U({\mathfrak g}) \otimes_{U(\fh)} 
\Lambda^i({\mathfrak g}/ \fh)$.  Then, for any ${\mathfrak g}$-module $M$, the sequence (with suitably defined  maps)  $\cdots \rightarrow
P_2 \otimes_{\C} M \rightarrow P_1 \otimes_{\C} M \rightarrow M \rightarrow 0$ is a $({\mathfrak g}, \fh)$-projective resolution \cite[3.1.8]{Kum:02}.  If $M=M(\lambda)$, then $P_i \otimes_{\C} M(\lambda) \in \W$, and so this gives a projective resolution in $\W$.  We then have
\begin{eqnarray*}
\Ext_{\W}^i \left( M(\lambda), L(\nu) \right) &\cong& \Ext_{({\mathfrak g}, \fh) }^i \left( M(\lambda), L(\nu) \right) \\
&\cong& \Ext_{(\fb^+, \fh)}^i \left(\C_\lambda, L(\nu) \right) \\
&\cong& \Ext_{(\fb^+, \fh)}^i \left( \C, \C_{-\lambda}\otimes L(\nu)\right) \\
&\cong& \operatorname{H}^i \left( \fb^+, \fh, \C_{-\lambda}\otimes L(\nu)\right) \\
&\cong& \operatorname{H}^i (\fn^+, L(\nu) )^{\lambda}.
\end{eqnarray*}
The second through fourth isomorphisms follow from \cite[3.1.14, 3.1.13, 3.9]{Kum:02}, respectively.  The last isomorphism follows from definitions.
\end{proof}
There are two functors $\eta : \W \rightarrow \W(\mu)$ and $\theta : \W \rightarrow \W(\mu)$ defined as follows.  For $M \in \W$, there is a unique minimal submodule $M' \subseteq M$ such that $M/M' \in \W(\mu)$.  Define $\eta M = M / M'$.  Note that $M'$ is generated, as a ${\mathfrak g}$-module, by $\bigoplus_{\lambda \not\leq \mu} M^{\lambda}$. Then, for any $N, M \in \W$, and any ${\mathfrak g}$-module homomorphism $f: M \rightarrow N$, $f(M') \subseteq N'$.  Therefore, $f$ induces a homomorphism from $\eta M$ to $\eta N$.  Define $\eta (f)$ to be this map.  For $M \in \W$, there is also a unique maximal submodule $M''$ such that $M'' \in \W(\mu)$.  Define $\theta M = M''$.  For any $N, M \in \W$, and any ${\mathfrak g}$-module homomorphism $f: M \rightarrow N$, define $\theta(g) = g \mid_{M''}$.  Using these functors we relate $\Ext^i_{\W}(-,-)$ and $\Ext_{\W(\mu)}^i(-,-)$.

\begin{lem}[B] \Label{L:WOmu}
Let $M, N \in \W(\mu)$.  Then $\Ext_{\W(\mu)}^i(M,N)= \Ext_{\W}^i(M,N)$.
\end{lem}
\begin{proof}
First observe that $\eta$ takes projectives to projectives and $\theta$ takes injectives to injectives. 
Let $N \rightarrow I_{\bullet}$ be an injective resolution in $\W$.  Since $M \in \W(\mu)$, $\Hom_{\fg} (M, I_k) \cong \Hom_{\fg} (M, \theta I_k)$.  Therefore, $\Ext^i_{\W}(M,N)= \mbox{H}^i(\Hom_{\fg}  (M, I_k)) \cong \mbox{H}^i(\Hom_{\fg} (M, \theta I_k))$. The lemma follows if we can show that $\theta $ is acyclic on $N$ because this would imply that 
$N=\theta N \rightarrow \theta I_{\bullet}$ is an injective resolution. 

Note that $\mbox{H}^i(\theta I_{\bullet})=0$ if and only if $\mbox{H}^i((\theta I_{\bullet})^{\gamma})=0$ for all $\gamma \in \fh^*$.  For $\gamma \in \fh^{*}$, define $P_\gamma = U({\mathfrak g}) \otimes_{U(\fh)} \C_{\gamma}$.  Then $(\theta I_k)^{\gamma}=\Hom_{\fg} (P_{\gamma}, \theta I_k)$.  This implies
\begin{eqnarray*}
\mbox{H}^i((\theta I_{\bullet})^{\gamma})&\cong&\mbox{H}^i(\Hom_{\fg} (P_{\gamma}, \theta I_\bullet)) \\
&\cong&\mbox{H}^i(\Hom_{\fg} (\eta P_{\gamma}, I_\bullet)) \\
&\cong&\Ext^i_{\W}(\eta P_{\gamma}, N).
\end{eqnarray*}
Therefore, to complete the proof it is enough to show that $\Ext^i_{\W}(\eta P_{\gamma}, N)=0$ for $i\geq 1$.  There is a short exact sequence $0 \rightarrow P_{\gamma}' \rightarrow P_{\gamma} \rightarrow \eta P_{\gamma} \rightarrow 0$, which give a long exact sequence
$$
\cdots \rightarrow \Ext^{i-1}_{\W} (P_{\gamma}', N)  \rightarrow \Ext^{i}_{\W} (\eta P_{\gamma}, N)  \rightarrow \Ext^{i}_{\W} (P_{\gamma}, N) \rightarrow \cdots
$$
Since $P_{\gamma}$ is projective in $\W$, $\Ext^{i}_{\W} (P_{\gamma}, N)=0$ for $i \geq 1$.  We claim $\Ext^{i-1}_{\W} (P_{\gamma}', N)=0$ for all $i$.  To see this, let $P^{\not\leq \mu}_{\gamma} = \mbox{span} \{m \in P_{\gamma} \mid m \in P_{\gamma}^{\nu} \ \mbox{for some $\nu \not\leq \mu$}, \ m=x \otimes 1 \ \mbox{for some $x \in U(\fn^+)$} \}$.   Then  $P^{\not\leq \mu}_{\gamma}$ is a $\fb^+$-module.   

Define $\W_{\fb^+}$ to be the category of $\fb^+$-modules $M$ such that $M= \bigoplus_{\lambda \in \fh^{*}} M^{\lambda}$.  Define $Q_k= U(\fb^+) \otimes_{U(\fh)} \Lambda^k (\fb^+/\fh) \otimes_{\C} P_{\gamma}^{\not\leq \mu}$.  Then $Q_{\bullet} \rightarrow P^{\not\leq \mu}_{\gamma}$  is a projective resolution of $P^{\not\leq \mu}_{\gamma}$ in $\W_{\fb^+}$ (cf.\ \cite[3.1.8]{Kum:02}).  Also, $U(\fg) \otimes_{U(\fb^+)} Q_{\bullet}$ is a projective resolution of $P_{\gamma}' \cong U(\fg) \otimes_{U(\fb^+)} P^{\not\leq \mu}_{\gamma}$ in $\W$.  Moreover,
$\Hom_{\fg} (U(\fg) \otimes_{U(\fb^+)} Q_{k}, N) = \Hom_{\fb^+} (Q_{k}, N)=0$ since $Q_k^{\nu} \neq 0$ only for $\nu \not\leq \mu$ and $N^{\nu} \neq 0$ only for $\nu \leq \mu$.  Therefore, $\Ext^{i-1}_{\W} (P_{\gamma}', N)=0$ for all $i$.  This implies that $\Ext^i_{\W}(\eta P_{\gamma}, N)=0$ for $i\geq 1$.
\end{proof}

We now transfer the information from $\W(\mu)$ to $\mathcal C$.

\begin{thm} \Label{T:ExtML} Let $\mathcal C$ be a finite block of $\mathcal O$, an infinite block of $\O$ with a maximal element, or a truncated infinite block with minimal element. Let $\Lambda$ be the weight poset of $\mathcal C$ with maximal element $\mu$, and let $\lambda, \nu \in \Lambda$.  Then for $i \geq 0$,
\begin{itemize}
\item[(a)] if $\mathcal C$ is a thick block or a finite block
\begin{equation} \Label{E:ExtML}
\Ext_{\mathcal C}^i \left( M(\lambda), L(\nu) \right)\cong
\begin{cases}
\C &\mbox{if} \ \la \preceq \nu,\ l(\lambda)-l(\nu)=i ;\\
0 &\text{otherwise}.
\end{cases}
\end{equation}
\item[(b)] if $\mathcal C$ is a thin block
\begin{equation}  \Label{E:ExtML-1}
\Ext_{\mathcal C}^i \left( M(\lambda), L(\nu) \right)\cong
\begin{cases}
\C &\mbox{if} \ \la \preceq \nu,\ l(\lambda)-l(\nu)=i,\ i\leq 1;\\
0 &\text{otherwise}.
\end{cases}
\end{equation}
\end{itemize}
\end{thm}
\begin{proof}
Let $M, N \in \mathcal C$.  Given Proposition \ref{T:n-homology} and Lemmas \ref{L:relativecohom}(A) and (B), it is enough to show that $\Ext^i_{\mathcal C} (M,N) \cong \Ext^i_{\W(\mu)}(M,N)$.

Let $\gamma \in \fh^{*}$, and recall $P_\gamma = U({\mathfrak g}) \otimes_{U(\fh)} \C_{\gamma}$.
Then $P_\gamma$ is projective in $\W$, and thus $\eta P_\gamma$ is projective in $\W(\mu)$.  Moreover, $\eta P_{\gamma}$ is finitely generated: $\eta P_{\gamma}$ is generated by $1 \otimes 1$ if $\gamma \leq \mu$ and $\eta P_{\gamma}=0$ otherwise.  Thus, $\eta P_{\gamma} \in \O$.  Therefore, we can construct a resolution  $P_{\bullet} \rightarrow M$ of $M$ so that $P_i = \bigoplus_{j=1}^{n_i} \eta P_{\gamma^i_j} \in \O$  for some $\gamma^i_j \in \fh^{*}$, which is projective in $\W(\mu)$.  

Recall that modules in $\O$ decompose according to blocks.  Let $\tilde{P}_i$ be the component of $P_i$ contained in $\mathcal C$.  (If $\mathcal C$ is a truncated block, the component of $P_i$ corresponding to the full block will be contained in the truncation $\mathcal C$ since $P_i \in \W(\mu)$ and $\mathcal C$ is truncated at $\mu$.) Then,  $ \tilde{P}_{\bullet} \rightarrow M$ is a projective resolution in $\W(\mu)$ and lies entirely in $\mathcal C$.  This proves the theorem.
\end{proof}

\subsection{\bf Cohomology and Quotient Categories} \Label{S:Quotient}
Throughout this section, let
${\mathcal C}$ be an infinite block for the Category ${\mathcal O}$ with a maximal element. Then ${\mathcal C}$ is a highest
weight category which contains enough projective objects. Let $\Lambda$
be the corresponding weight poset indexing the simple objects
in ${\mathcal C}$. For $\lambda\in \Lambda$ let $P(\lambda)$
be the projective cover of $L(\lambda)$. Set $P=\oplus_{\lambda\in \Lambda} P(\lambda)$.
Then $P$ is a progenerator for ${\mathcal C}$ and ${\mathcal C}$ is
Morita equivalent to $\text{Mod}(B)$ where $B=\text{End}_{\mathcal C}(P)^{\op}$.

We now apply results as described in \cite[Thm.\ 3.5]{CPS1}. Let $\Omega$ be a finite coideal, that is, $\Omega=\Lambda- \{\gamma \in \Lambda \mid \gamma \preceq \delta \}$ for some fixed $\delta \in \Lambda$. 
Consider $P_{\Omega}=\oplus_{\lambda\in \Omega}P(\lambda)$ and
set $A=\text{End}_{\mathcal C}(P_{\Omega})^{\op}$. Then
there exists an idempotent $e\in B$, corresponding to the
sum of identity maps in $\text{End}_{\mathcal C}(P(\lambda))$ with
$\lambda\in \Omega$, such that $eBe=A$. Also, observe that
the quotient category ${\mathcal C}(\Omega)=\text{Mod}(A)$ is a highest weight
category.

For $\lambda\in \Omega$, set
$M_{\Omega}(\lambda)=eM(\lambda)$, $L_{\Omega}(\lambda)=eL(\lambda)$, $P_{\Omega}(\lambda)=eP(\lambda)$. Note
that $P_{\Omega}(\lambda)$ is the projective cover of $L_{\Omega}(\lambda)$.
The following proposition compares extensions between
Verma modules and simple modules in ${\mathcal C}$ and ${\mathcal C}(\Omega)$.
We remark that this result appears as \cite[Cor.\ 3.5]{CPS6} with more
finiteness restrictions.

\begin{prop}[A] \Label{P:VermaSimpleExtCompare} Let $\lambda,\nu \in \Omega$.
For all $i\geq 0$, 
$$\operatorname{Ext}_{\mathcal C}^{i}(M(\lambda),L(\nu))\cong
\operatorname{Ext}_{{\mathcal C}(\Omega)}^{i}(M_{\Omega}(\lambda),L_{\Omega}(\nu)).$$
\end{prop}

\begin{proof} Let $\lambda,\nu\in \Omega$. According to \cite[Thm.\ 2.2]{DEN}, there exists a first quadrant
spectral sequence,
\begin{equation*}
E_{2}^{i,j}=\text{Ext}^{i}_{B}(\text{Tor}^{A}_{j}(Be,M_{\Omega}(\lambda)),L(\nu))
\Rightarrow \text{Ext}^{i+j}_{A}
(M_{\Omega}(\lambda),L_{\Omega}(\nu)).
\end{equation*}
By \cite[Thm.\ 4.5]{DEN}, $\text{Tor}_{0}^{A}(Be,M_{\Omega}(\lambda))=M(\lambda)$.
We need to show that $\text{Tor}_{j}^{A}(Be,M_{\Omega}(\lambda))=0$ for $j\geq 1$.
Then the spectral sequence above collapses and yields
$$\text{Ext}^{i}_{B}(M(\lambda),L(\nu))\cong \text{Ext}^{i}_{A}(M_{\Omega}(\lambda),L_{\Omega}(\nu)).$$
for $i\geq 0$ and $\lambda,\nu\in \Omega$, as required.

First we consider the case when $j=1$. Since ${\mathcal C}(\Omega)$ is
a highest weight category, we may invoke \cite[Thm.\ 4.5]{DEN} which
states that $\text{Tor}_{0}^{A}(Be,M_{\Omega}(\lambda))=M(\lambda)$
and $M(\lambda)$ belongs to ${\mathcal X}$ (cf.\ \cite[\S 3.1]{DEN} for a definition
of ${\mathcal X}$). Therefore, $\text{Tor}_{1}^{A}(Be,M_{\Omega}(\lambda))=0$ by
\cite[Prop.\ 3.1(A)]{DEN}.

We now use induction on the ordering on the weights in $\Omega$ to show that
$\text{Tor}_{j}^{A}(Be,M_{\Omega}(\lambda))=0$ for $j\geq 2$.
If $\lambda$ is a maximal weight (relative to $\preceq$, the ordering introduced in \eqref{E:Ordering}) then $M_{\Omega}(\lambda)$ is the
projective cover of $L_{\Omega}(\lambda)$ and
$\text{Tor}_{j}^{A}(Be,M_{\Omega}(\lambda))=0$ for $j\geq 1$.
Now suppose that $\text{Tor}_{j}^{A}(Be,M_{\Omega}(\mu))=0$ for $j\geq 1$
for all $\mu\succ\lambda$, $\mu\in\Omega$. Consider the short exact sequence
$$0\rightarrow N\rightarrow P_{\Omega}(\lambda)\rightarrow M_{\Omega}(\lambda)
\rightarrow 0,$$
where $N$ has a filtration by modules $M_{\Omega}(\mu)$ with $\mu\succ\la$.
This induces a long exact sequence
$$\dots \leftarrow \text{Tor}_{j-1}^{A}(Be,N)\leftarrow
\text{Tor}_{j}^{A}(Be,M_{\Omega}(\lambda))\leftarrow
\text{Tor}_{j}^{A}(Be,P_{\Omega}(\lambda))\leftarrow \dots
$$
For $j\geq 1$, $\text{Tor}_{j}^{A}(Be,P_{\Omega}(\lambda))=0$, and for
$j\geq 2$, $\text{Tor}_{j-1}^{A}(Be,N)=0$ by the induction hypothesis. Thus from the long exact sequence
we can conclude for $j\geq 2$, $\text{Tor}_{j}^{A}(Be,M_{\Omega}(\lambda))=0$.
\end{proof}

Let $L(\lambda)$ and $L(\nu)$ be simple $B$-modules with $\lambda,\nu
\in \Omega$. Then $eL(\lambda)\neq 0$ and $eL(\nu)\neq 0$. Let
$$\dots \rightarrow P_{2}\rightarrow P_{1} \rightarrow P_{0} \rightarrow L(\lambda) \rightarrow 0$$
be the minimal projective resolution of $L(\lambda)$ in ${\mathcal C}$.
Set $\Omega^{n+1}(L(\lambda))$ to be the kernel of the map
$P_{n}\rightarrow P_{n-1}$. By convention, we let $\Omega^{0}(L(\lambda))=L(\lambda)$.
Under a suitable condition on the
minimal projective resolution, we can compare extensions between
simple modules in ${\mathcal C}$ and ${\mathcal C}(\Omega)$.  This comparison depends on bounding the composition factors in the projective resolution of $L(\la)$.  The following proposition provides such a bound.

\begin{prop}[B] \Label{P:radMLextbound} Let $\la, \nu \in \Lambda$.
\begin{itemize}
\item[(i)] If $\Ext^n_{\CC}(\rad M(\la),L(\nu)) \ne 0$ then $l(\la)-l(\nu)\le n-1$.
\item[(ii)] If $\Ext^n_{\CC}(L(\la),L(\nu)) \ne 0$ then $l(\la)-l(\nu)\le n$.
\end{itemize}
\end{prop}
\begin{proof}
In this proof we assume that $\mathcal C$ is a thick block.  In the case that $\mathcal C$ is a thin block or a finite block, the proposition follows from similar arguments.

(i) We prove this by induction on $n$. Let $n=0$. Since $\rad M(\la) = M(\la_1)+M(\la'_1)$ (notation as in Figure~\ref{embedd4}), $\Hom_{\CC}(\rad M(\la), \linebreak[0] L(\nu)) \ne 0$ if and only if $\nu=\la_1$ or $\la'_1$, whence $l(\la)-l(\nu)=-1$.

Assume the result is true for $n-1$ and all pairs of weights in $\Lambda$. Recall from \eqref{E:VMrad} that $\rad M(\la)=M(\la_1)+M(\la'_1)$ and $\rad^2 M(\la)=M(\la_2)+M(\la'_2)=\rad M(\la_1)=\rad M(\la'_1)$.  Thus we have a short exact sequence
$$
0 \to \rad M(\la_1) = \rad M(\la'_1) \to M(\la_1)\oplus M(\la'_1) \to \rad M(\la) \to 0,
$$
where the inclusion sends $x$ to $(x,-x)$ and the surjection sends $(x,y)$ to $x+y$.
This induces a long exact sequence
$$
\cdots \to \Ext^{n-1}_{\CC}(\rad M(\la_1),L(\nu)) \to
\Ext^n_{\CC}(\rad M(\la),L(\nu)) \to \Ext^n_{\CC}(M(\la_1)\oplus M(\la'_1),L(\nu)) \to \cdots.
$$
Suppose $l(\la)-l(\nu)>n-1$. Then $l(\la_1)-l(\nu)=l(\la)+1-l(\nu)>n$.  This implies
$\Ext^n_{\CC}(M(\la_1),L(\nu))=0$ by Theorem \ref{T:ExtML}, and similarly for $\la'_1$. Also, $l(\la_1)-l(\nu)>(n-1)+1$, so
$\Ext^{n-1}_{\CC}(\rad M(\la_1),\linebreak[0]L(\nu))=0$ by induction. This implies $\Ext^n_{\CC}(\rad M(\la),L(\nu))=0$.

(ii) The proof is again by induction on $n$. The result is clear for $n=0$. Assume it is true for $n-1$.
Consider the short exact sequence
$$
0 \to \rad M(\la) \to M(\la) \to L(\la) \to 0.
$$
This induces a long exact sequence
$$
\cdots \to \Ext^{n-1}_{\CC}(\rad M(\la),L(\nu)) \to \Ext^n_{\CC}(L(\la),L(\nu)) \to
\Ext^n_{\CC}(M(\la),L(\nu)) \to \cdots.
$$
Suppose $l(\la)-l(\nu)>n$. Then $\Ext^{n-1}_{\CC}(\rad M(\la),L(\nu))=0$ by part (i),
and $\Ext^n_{\CC}(M(\la),\linebreak[0] L(\mu))=0$ by Theorem \ref{T:ExtML}. This implies $\Ext^n_{\CC}
(L(\la),L(\nu))=0$.
\end{proof}
Recall $P_{\bullet} \rightarrow L(\la)$ is a minimal projective resolution of $L(\la)$.  For $\gamma \in \Lambda$, if $L(\gamma)$ is a composition factor of $\hd P_{j}$, then $\Ext^j_{\CC}(L(\lambda), L(\gamma)) \neq 0$.  Therefore, Proposition \ref{P:radMLextbound}(B) gives a bound on the composition factors which can appear in $\hd P_{j}$.  This is the condition needed to compare extensions between
simple modules in ${\mathcal C}$ and ${\mathcal C}(\Omega)$.

\begin{prop}[C] \Label{P:SimpleExtCompare} 
Let $\la, \nu \in \Omega$ and  define $N= \mbox{min} \{\, | l(\la)-l(\gamma)| \ \: \gamma \in \Lambda - \Omega \,\} -1$.  Then, for $j=0,1, \ldots, N$,
$$\operatorname{Ext}^{j}_{{\mathcal C}}(L(\lambda),L(\nu))
\cong \operatorname{Ext}^{j}_{{\mathcal C}(\Omega)}(L_{\Omega}(\lambda),L_{\Omega}(\nu)).
$$
\end{prop}

\begin{proof} We first claim that $BeP_{j}=P_{j}$ for $j=0, \ldots, N$.  
Note that $BeP_{j}=P_{j}$ if and only if $\hd P_{j}$ contains no composition factors which
are killed by the idempotent $e$.   
Suppose that $L(\gamma) \subseteq \hd P_{j}$.  Then $\Ext_{\mathcal C}^j (L(\gamma), L(\lambda))  \neq 0$.  From the proof of Theorem \ref{T:ExtML} and using the duality on $\W(\mu)$, we have that 
$$
\Ext_{\mathcal C}^j (L(\gamma), L(\lambda)) \cong \Ext_{\W(\mu)}^j (L(\gamma), L(\lambda)) \cong \Ext_{\W(\mu)}^j (L(\lambda), L(\gamma)) \cong  \Ext_{\mathcal C}^j (L(\lambda), L(\gamma)).
$$
Then Proposition \ref{P:radMLextbound}(B) implies that $| l(\lambda)-l(\gamma) | \leq j$.  Therefore, for $j \leq N$, $\gamma \in \Omega$, and so $eL(\gamma) \neq 0$. Since
$\hd \Omega^{j}(L(\lambda))
\cong \hd P_{j}$,
we have $Be\Omega^{j}(L(\lambda))=\Omega^{j}(L(\lambda))$ for $j=0,1,2,\dots,N$.

Let $j=0,1,\dots,N-1$. Since $Be\Omega^{j}(L(\lambda))=\Omega^{j}(L(\lambda))$
there exists a short exact sequence \cite[Thm.\ 3.2]{DEN},
\begin{equation*}
0\rightarrow \Omega^{j+1}(L(\lambda))/Be\Omega^{j+1}(L(\lambda))
\rightarrow \operatorname{Tor}_{0}^{A}(Be,e\Omega^{j}(L(\lambda)))\rightarrow
\Omega^{j}(L(\lambda))\rightarrow 0.
\end{equation*}
Note that we are using the fact that $\Omega^{1}(\Omega^{j}(L(\lambda)))\cong \Omega^{j+1}(L(\lambda))$.
Since $\Omega^{j+1}(L(\lambda))=Be\Omega^{j+1}(L(\lambda))$, we have
\begin{equation*}
\operatorname{Tor}_{0}^{A}(Be,e\Omega^{j}(L(\lambda)))\cong \Omega^{j}(L(\lambda)).
\end{equation*}
Finally, we can apply \cite[Thm.\ 2.4(B)(ii)]{DEN} and a dimension shifting argument
twice (cf.\ \cite[Cor.\ 2.5.4]{Ben:98a}) to see that
\begin{eqnarray*}
\text{Ext}^{j+1}_{B}(L(\lambda),L(\nu))&\cong& \text{Ext}^{1}_{B}(\Omega^{j}(L(\lambda)),L(\nu))\\
&\cong& \text{Ext}^{1}_{B}(\operatorname{Tor}_{0}^{A}(Be,e\Omega^{j}(L(\lambda))),L(\nu))\\
&\cong& \text{Ext}^{1}_{A}(e\Omega^{j}(L(\lambda)),eL(\nu))\\
&\cong& \text{Ext}^{1}_{A}(\Omega^{j}(eL(\lambda)),eL(\nu))\\
&\cong& \text{Ext}^{j+1}_{A}(eL(\lambda),eL(\nu))\\
&\cong& \text{Ext}^{j+1}_{A}(L_{\Omega}(\lambda),L_{\Omega}(\nu))
\end{eqnarray*}
In the identifications above to justify the step between lines 3 and 4,
observe that the idempotent functor $e(-):\text{Mod}(B)\rightarrow
\text{Mod}(A)$ is exact. Moreover, $BeP_{j}=P_{j}$ for $j=1,2,\dots, N$ so we have an
exact sequence of projective $A$-modules:
$$eP_{N}\rightarrow eP_{N-1}\rightarrow \dots eP_{1} \rightarrow eP_{0} \rightarrow eL(\lambda) \rightarrow 0.$$
This implies that $e\Omega^{j}(L(\lambda))\cong \Omega^{j}(eL(\lambda))\oplus Q_{j}$ where $Q_{j}$ is a projective $A$-module for $j=0,1,\dots,N$.
\end{proof}

\subsection{\bf Extensions Between Simple Modules}
Let $\mathcal C$ be a finite block, an infinite block with a maximal element, or a truncation of an infinite block with a minimal element.   Let $\Lambda$ be the weight poset of $\mathcal C$, with length function $l: \Lambda \rightarrow \Z$.  
\begin{thm} \Label{T:ExtLL}
Let $\lambda, \nu \in \Lambda$.  Then,
\begin{itemize}
\item[(a)] if $\mathcal C$ is a thick block or a finite block
\begin{equation}
 \dim \Ext^n_{\mathcal C} (L(\lambda), L(\nu)) = \# \{ \gamma \in \Lambda \mid \gamma \preceq \la, \nu;\ 2l(\gamma)-l(\lambda)-l(\nu)=n \}
\end{equation}
\item[(b)] if $\mathcal C$ is a thin block
\begin{equation}
 \dim \Ext^n_{\mathcal C} (L(\lambda), L(\nu)) = 
 \# \{ \gamma \in \Lambda \mid \gamma \preceq \la, \nu;\ 2l(\gamma)-l(\lambda)-l(\nu)=n \}
\end{equation}
if $n \leq 2$, and equals zero otherwise.
\end{itemize}
In particular, $ \Ext^n_{\mathcal C} (L(\lambda), L(\nu)) \neq 0$ only when $n \equiv \left( l(\lambda) - l(\nu) \right) \pmod{2}$.
\end{thm}
\begin{proof}
Suppose $\mathcal C$ is an infinite block with a maximal element.  Let $\Omega$ be a finite coideal of $\Lambda$ containing $\la, \nu$.  From Proposition \ref{P:VermaSimpleExtCompare}(A), we know $\Ext^i_{\mathcal C} (M(\lambda), L(\nu))= \Ext^i_{{\mathcal C}(\Omega)} (M_{\Omega}(\lambda),\linebreak[0] L_{\Omega}(\nu))$.  For a fixed $n \in \Z_{>0}$, we assume that $\Omega$ is sufficiently large so that $\gamma \in \Omega$ for all $\gamma \in \Lambda$ with  $| l(\la) - l(\gamma) | \leq n$.  Then Proposition \ref{P:SimpleExtCompare}(C) implies that $\Ext^n_{\mathcal C} (L(\lambda), L(\nu))
\cong \Ext^n_{{\mathcal C}(\Omega)} (L_{\Omega}(\lambda), L_{\Omega}(\nu))$.

Thus, by replacing $\mathcal C$ by a quotient category ${\mathcal C}(\Omega)$ where appropriate, we may assume that $\mathcal C$ is a highest weight category with finite weight poset $\Lambda$.  Because the objects of $\mathcal C$ have finite composition length, $\mathcal C$ is closed under the duality $D$ on $\W(\mu)$.  Define $A(\gamma)=DM(\gamma)$.

Now apply \cite[3.5]{CPS3}:
\begin{equation} \Label{E:PF}
 \dim \Ext^n_{\mathcal C} (L(\lambda), L(\nu))= \sum_{\gamma \in \Lambda, i, j \in \Z_{\geq 0}, i+j=n}  \dim \Ext^i_{\mathcal C} (L( \lambda), A(\gamma)) \dim \Ext^j_{\mathcal C} (M(\gamma), L(\nu)).
\end{equation}
  Using the duality on $\mathcal C$, $\Ext^i_{\mathcal C} (L(\la), A(\gamma))\cong \Ext^i_{\mathcal C} (M(\gamma), L(\la))$.  Then Theorem \ref{T:ExtML} gives the result.
\end{proof}

\subsection{\bf $\Ext^{1}$-quivers} \Label{SS:Ext1quiver}
Let $\mathcal C$ be a finite block, a quotient of an infinite block with a maximal element, or a truncation of an infinite block with a minimal element.   Let $\Lambda$ be the (finite) weight poset of $\mathcal C$. The $\Ext^{1}$-quiver of $\mathcal C$ is defined to be the directed graph with vertices labelled by $\Lambda$, and with $\dim \Ext^{1}_{\mathcal C}(L(\lambda),L(\mu))$ edges from $\lambda$ to $\mu$. It is clear from Theorem \ref{T:ExtLL} and Proposition \ref{P:SimpleExtCompare} (C) that the $\Ext^{1}$-quiver of $\mathcal C$ is obtained from the poset $\Lambda$ simply by replacing each edge by a pair of directed edges, one pointing in each direction.

The edges from $\lambda$ to $\mu$ can also be viewed as representing linearly independent elements of $\Hom_{\mathcal C}(P(\lambda),P(\mu))$ in the finite dimensional algebra 
$$
A=\End_{\mathcal C}\left(\bigoplus_{\la\in\Lambda} P(\la)\right)^{\text{op}}.
$$
One can ask for the relations that exist between the maps in this algebra, which provides a presentation 
of the algebra by the quiver with relations.

Suppose that $\mathcal C$ is either a finite block, or a finite quotient or truncation of a thin block. Then $\Lambda$ is a simple chain, say of length $n$, and it is quite easy to write down the structure of the projective indecomposables $P(\la)$. This is done for $n=4$ in \cite{FNP}, and the pattern is the same for any $n$. Moreover, if the elements of $\Lambda$ are numbered $\la_{1},\dots,\la_{n}$ from top to bottom, and if $\a_{i}$ (resp.\ $\b_{i}$) represents the map from $P(\la_{i})$ to $P(\la_{i+1})$ (resp.\ $P(\la_{i+1})$ to $P(\la_{i})$) for $1\le i\le n-1$, then one sees easily that the following relations hold (up to scalar multiples):
\begin{equation} \Label{E:thinquiverrelations}
\a_{1}\b_{1}=0, \qquad \b_{i}\a_{i}=\a_{i+1}\b_{i+1} \text{ for } 1\le i \le n-2
\end{equation}
(note that maps compose left-to-right, because of the $( )^{\text{op}}$ in the definition of $A$). 

Now we can assume we're working in the basic algebra with simple modules (resp., projective indecomposable modules) 
labelled by $\widehat{L}(\lambda_{i})$ (resp., $\widehat{P}(\lambda_{i})$), $i=1,2,\dots, n$. Note that 
$\dim \widehat{L}(\la_{i})=1$ for every $i$ so that $\dim A=\sum_{i=1}^{n} \dim \widehat{P}(\la_{i})$, which is easy to 
compute using the known structures of the $\widehat{P}(\la_{i})$. On the other hand, using the relations given in 
\eqref{E:thinquiverrelations}, one can check directly that there are at most $\sum_{i=1}^{n} \dim \widehat{P}(\la_{i})$ linearly 
independent words in the $\a_{i}$ and $\b_{i}$. Thus \eqref{E:thinquiverrelations} must be {\it all} the relations. 

In contrast, suppose that $\mathcal C$ is a finite quotient or truncation of a thick block. Then the poset $\Lambda$ is isomorphic to the Bruhat order on a dihedral group. In this case the structure of the projectives, and the exact nature of the relations, seem to be quite difficult to deduce. For example, Stroppel in \cite{Str} works out the relations for the $\Ext^{1}$-quiver of the regular blocks of Category $\mathcal O$ for the finite simple complex rank 2 Lie algebras, using some deep results of Soergel.  Not only are the answers quite complicated (e.g., for $G_{2}$ there are 70 relations), but as far as we are aware the analogs of Soergel's results are not known for the Virasoro algebra. Nonetheless, based on Stroppel's computations, we speculate that the relations in the $\Ext^{1}$-quiver of $\mathcal C$ are all quadratic.

\section{Kazhdan-Lusztig Theories and Koszulity}
\Label{S:KL}

\subsection{} Let $B= B_0 \oplus B_1 \oplus \cdots \oplus B_{q}$ be a finite-dimensional 
graded algebra, and let ${\mathcal C}^{gr}_B$ 
be the category of finite-dimensional graded $B$-modules. Regard every simple $B$-module 
$L$ as concentrated in degree zero; then the simple modules in ${\mathcal C}^{gr}_B$ can 
be obtained by shifting the gradings of the simple $B$-modules. If $L$ is a simple $B$-module 
then $L(i)$ will denote the simple module in ${\mathcal C}^{gr}_B$ by shifting $i$ places 
to the right (cf.\ \cite[\S 1.3]{CPS5}). The algebra $B$ is {\it Koszul} if for all 
simple $B$-modules $L$ and $L'$ and $m,n,p \in \Z$,
\begin{equation}
\operatorname{Ext}^p_{{\mathcal C}^{gr}_B} (L(m), L'(n)) \neq 0 \Rightarrow n-m=p.
\end{equation} 

Now let $\mathcal C$ be either a finite block of $\O$, a truncation of an infinite block of $\O$ with 
a minimal element, or a quotient of an infinite block of $\O$ with a maximal element.  Then $\mathcal C$ 
is a highest weight category (with duality) having a finite weight poset $\Lambda$ and length function $l$.  
Moreover, Proposition  \ref{T:ExtML} implies that $\mathcal C$ has a Kazhdan-Luzstig theory, as defined in 
\cite{CPS3}. Recall that $\mathcal C$ is equivalent to $\mbox{Mod}(A)$ for a finite-dimensional algebra $A$.  
Let $\gr A$ be the associated graded algebra obtained by using the radical filtration on $A$. 
Moreover, set $L = \bigoplus_{\lambda \in \Lambda} L(\lambda)$, and define  the homological dual of $A$ to be 
$A^!= \Ext^{\bullet}_{\mathcal C} (L, L)$. The following theorem establishes Koszulity results on $A$. 

\begin{thm} Let $\mathcal C$ be as described above and $A$ be the associated quasi-hereditary algebra. Then 
\begin{itemize} 
\item[(a)] $A^{!}$ is Koszul
\item[(b)] $\operatorname{gr }A$ is Koszul
\item[(c)] $(A^{!})^{!}\cong \operatorname{gr }A$
\end{itemize} 
\end{thm} 

\begin{proof} In order to prove the theorem, it suffices to check the condition
\begin{equation}
\Ext^n_{\CC}(\rad^i M(\la),L(\mu))\ne 0 \Rightarrow n \equiv l(\mu)-l(\la)+i\!\!\!\!\pmod 2\ \text{ for all } 
\la,\mu\in\Lambda.\tag{SKL${}'$}
\end{equation}
In principle we should also check the same parity vanishing for $\Ext^n_{\CC}(L(\la),A(\mu)/\soc^i A(\mu))$ 
but this follows by duality in our setting.

Once ($\text{SKL}^{\prime}$) is established then by \cite[Lemma 2.1.5]{CPS5} ${\mathcal C}^{gr}_{A^{!}}$ 
has a graded KL-theory, so ${\mathcal C}^{gr}_{A^{!}}$ and $A^{!}$ are Koszul using \cite[Thm.\ 3.9]{CPS4}. 
The condition ($\text{SKL}^{\prime}$) implies parts (b) and (c) by \cite[Thm.\ 2.2.1]{CPS5}. 

The case $i=0$ of (SKL${}'$) follows immediately from Theorem \ref{T:ExtML}. Assume $i>0$. Then $\rad^i M(\la)$ 
is either $0$, a Verma module $M(\nu)$ with $l(\nu)-l(\la)=i$, or a sum of two such Verma modules. The first case 
is trivial, and in the second, we can use the same argument as for $i=0$. So assume we are in the third case. 
Let $\la_i$ and $\la'_i$ be the two elements $\nu$ satisfying $l(\nu)-l(\la)=i$. We have a short exact sequence
$$
0\to M(\la_i) \to \rad^iM(\la) \to L(\la'_i) \to 0.
$$
The corresponding long exact sequence is
$$
\cdots \to \Ext^n_{\CC}(L(\la'_i),L(\mu)) \to \Ext^n_{\CC}(\rad^i M(\la),L(\mu)) \to \Ext^n_{\CC}(M(\la_i),L(\mu)) 
\to \cdots
$$
and we are assuming the middle term is nonzero. Then one of the two adjacent terms must be nonzero. If the term 
on the right is nonzero, then the same argument as for $i=0$ gives the desired parity condition, since $l(\la_i)
=l(\la)+i$. If the term on the left is nonzero, then by Theorem \ref{T:ExtLL} we have $n\equiv l(\la'_i)-l(\mu)
\equiv l(\la)+i-l(\mu)\pmod 2$, which is equivalent to the required condition.
\end{proof}

\begin{rems}
1. Since this proof holds for all quotient categories of $\mathcal C$, this shows that $\mathcal C$ 
satisfies the strong Kazhdan Lusztig condition (cf.\ \cite[2.4.1]{CPS5}).

2. Suppose $\mathcal C$ comes from either a finite or a thin block of $\mathcal O$. Since the relations in the $\Ext^{1}$-quiver 
of $\mathcal C$ are all homogeneous (in fact quadratic (cf.\ \eqref{E:thinquiverrelations})), it follows that $A$ itself is tightly 
graded (i.e., $A\cong \operatorname{gr} A$). In particular in this case $A$ is Koszul.

3. An interesting open question 
is to determine whether $A$ is tightly graded or whether $A$ itself is Koszul when $A$ is associated to a thick block. The answers would be affirmative if the relations are all quadratic, as speculated in Section \ref{SS:Ext1quiver}.
\end{rems}

\let\section=\oldsection

\def\scr{\mathcal}\def\cprime{$'$} \def\germ{\mathfrak}

\end{document}